\documentclass[12
pt,twoside,a4paper]{amsart}
\usepackage{amssymb}
\usepackage[latin1]{inputenc}
\usepackage[T1]{fontenc}
\usepackage{times}

\def\ephi{e^{-\phi}}

\def\hpq0{h^{p,q}_{\leq 0}}
\def\Hpq0{\H_{\leq 0}^{p,q}}
\def\gr{\otimes}

\def\dbar{\bar\partial}
\def\ddbar{\partial\dbar}

\def\R{{\mathbb R}}

\def\C{{\mathbb C}}
\def\Cn{\C^n}

\def\K{{\mathcal K}}

\def\A{{\mathcal A}}

\def\H{{\mathcal H}}

\def\F{{\mathcal F}}

\def\Re{{\rm Re\,  }}

\def\L{{\mathcal L}}

\def\be{\begin{equation}}
\def\ee{\end{equation}}

\newtheorem{thm}{Theorem}[section]
\newtheorem{lma}[thm]{Lemma}
\newtheorem{cor}[thm]{Corollary}
\newtheorem{prop}[thm]{Proposition}

\theoremstyle{definition}

\theoremstyle{remark}

\newtheorem{preremark}{Remark}
\newtheorem{preex}{Example}

\numberwithin{equation}{section}

\begin{document}

\title[]
{Positivity of direct image bundles and convexity on the space of
Kähler metrics.}

\author[]{ Bo Berndtsson}

\address{B Berndtsson :Department of Mathematics\\Chalmers University
  of Technology 
  and the University of G\"oteborg\\S-412 96 G\"OTEBORG\\SWEDEN,\\}

\email{ bob@math.chalmers.se}

\begin{abstract}
{We  develop some results from \cite{B1} on the positivity of
  direct image bundles in the particular case of a trivial fibration
  over a one-dimensional base. We also apply the results to study
  variations of Kähler metrics. }
\end{abstract}

\bigskip

\maketitle

\section{Introduction}

In a previous paper, \cite{B1}, we have studied curvature properties
of vector bundles that arise as direct images of line bundles over a
Kähler fibration. Here we will continue this study in a very special case
- trivial fibrations over a one-dimensional base- and elaborate on its
connection with problems of variations of Kähler 
metrics on a compact manifold.

Let $Z$ be a compact complex $n$-dimensional manifold, and let $\hat
L$ be a positive 
line bundle over $Z$. Put $X=U\times Z$, where $U$ is a domain in
$\C$, and denote by $L$ the pullback of $\hat L$ under the projection
from $X$ to $Z$. Let $\phi$ be a hermitian metric on $L$. Then $\phi$
can be written
$$
\phi =\phi_0 +\psi,
$$
where $\phi_0=\phi_0(z)$ is some fixed metric on $\hat L$, pulled back to $X$,
and $\psi=\psi(t,z)$ is a function of $t$ in $U$ and $z$ in $Z$. We
can thus think of $\phi$ as a family of metrics on $\hat L$, indexed
by $t$ in $U$, or as a (complex) curve in the affine space of all
metrics on $\hat L$. Put
$$
\hat E= H^0(Z,\hat L\gr K_Z),
$$
the space of global holomorphic $\hat L$-valued $(n,0)$-forms on
$Z$. Denote by $E$ the (trivial) vector bundle over $U$ 
with fiber $\hat E$.

Even though $E$ is a trivial bundle it has a naturally defined metric,
$H_\phi $ which is not trivial. $H_\phi$ is defined by
$$
\|u\|^2_t = \int_{\pi^{-1}(t)} [u,u]e^{-\phi}=\int_Z [u,u]e^{-\phi_0
  -\psi(t,z)},
$$
if $u$ is an element in $E_t=\hat E$, the fiber of $E$ over $t$. Here
$
[u,u]e^{-\phi}$ is defined by
$$
[u,u]e^{-\phi}=c_n u_j\wedge\bar u_j e^{-\phi_j}
$$
if $u_j$ and $\phi_j$ are local representatives of $u$ and $\phi$ with
respect to a local trivialization. The constant $c_n=i^{n^2}$ is 
chosen to make this expression nonnegative. 

With this metric $E$ becomes an Hermitian vector bundle and it makes
sense to consider its curvature, $\Theta^E$. It is proved in \cite{B1}
that if $\phi$ is a (semi)positive metric on $L$, then $H_\phi$ is a
(semi)positive metric on $E$. Moreover, assuming that $Z$ has no
nontrivial global holomorphic vector fields, $\Theta^E$ is strictly
positive at a point $t$ unless  $\omega^t=i\ddbar_z\phi$ is stationary
so that
$$
\frac{\partial}{\partial t} \omega^t|_t =0.
$$

\bigskip

We are now ready to describe the content of this paper. In the next
section we will first reprove the formula for the curvature (Theorem
2.1) from \cite
{B1}, which is particularily simple in this situation. We then
proceed to give exact conditions on $\phi$ that
characterize when the curvature $\Theta^E$ has a null-vector, for the
case when $\phi$ itself is semipositive. This result is implicit in
\cite{B1}, but not explicitly stated there. Then we go on to show that if
the curvature is degenerate not only for one single $t$ but for all
$t$ in an open neighbourhood of 0, then
$$
\omega^t =T_{t}^*{\omega^{0}},
$$
where $T_t$ is the flow of some holomorphic vector field on $Z$.

In section 3 we then give a result corresponding to Theorem 2.1 for
a trivial vector bundle $F$ with fiber $ H^0(Z, L)$ instead of 
$H^0(Z,L\gr K_Z)$. In this case we define our $L^2$-norms on the fiber
by
$$
\|u\|^2=\int_Z |u|^2 \ephi\omega_n
$$
(see below). Surprisingly,  we then get an estimate from {\it above} of the
curvature, which implies among other things that if $\omega^{n+1}=0$ on
the total space $U\times Z$, then $\Theta^F$ is {\it seminegative}.

Next we consider the Hermitian bundles $E(p)$, defined in the same way
as $E$, but
replacing $L$ by $L^p$, for some positive integer $p$. We then take
$p\phi$ for our metric on $L^p$, and put
$\Theta^p=\Theta^{E(p)}$ for ease of notation. It is an a priori non
trivial fact, but an immediate
consequence of the characterization mentioned above, that if, at some
$t$ and for some $p_0$, $\Theta^{p_0}$ is degenerate, then $\Theta^p$
actually vanishes at $t$ for all $p$. 

Motivated in part by this observation we  then study the asymptotic
behaviour of $\Theta^p$ as $p$ tends to infinity. Let 
$$
\tau(p)= \mathbf{tr} \Theta^p/d_p,
$$
where $d_p$ is the rank of $E(p)$, i e the dimension  of $H^0(Z,\hat
L^p\gr K_Z)$. Thus $\tau(p)$ is the average of the eigenvalues of the
curvature. 
In section 4 we give an asymptotic formula for $\tau(p)$ containing
one term of order $p$, one term of  order zero and one term that
vanishes in the limit.
It follows from this formula that if we assume $i\ddbar\phi$ to be
semipositive, then
$$
\liminf_{p\rightarrow \infty} \tau(p)=0,
$$ 
implies that  the condition characterizing  degeneracy follows.
Hence, in this case,  $\Theta^p$ actually vanishes for all $p$.  Ma
and Zhang \cite{Ma-Zhang} have announced the existence of a 
complete asymptotic expansion for the (or 'a')  kernel of the full
curvature operator 
in a situation more general than
ours and also computed the first two terms on the diagonal.   I do not know if
our result (which is much more elementary) follows from theirs. Our
result does however follow from the relation between our vector
bundles $E(p)$ and the {\it Mabuchi functional} from Kähler geometry
that we  discuss in section 5. Hence the reader who so wishes can skip
directly to section 5 after a brief look at our recollection of the
Tian-Catlin-Zelditch asymptotic formula for the Bergman kernel in
section 4. 

In section 5
we relate our  results to the study of
variations of Kähler metrics on $Z$, cf \cite{Mabuchi},
\cite{Bando-Mabuchi}, \cite{Donaldson 1}, \cite{Donaldson 2} \cite{Phong-Sturm}
\cite{Chen}. In these works one considers 
the space of all Kähler metrics $\omega$ on $Z$ that lie in one fixed
cohomology class, viz, $c(\hat L)$,  the Chern class of $\hat L$. Any
such metric can be 
written
$$
\omega =\omega^0 +i\ddbar\psi,
$$
where $\omega^0$ is the fixed reference metric and $\psi$ is a
function on $Z$. The space of all such Kähler metrics
can therefore be identified with the space  $\mathcal{K}$ of all - say smooth -
metrics on $\hat L$, modulo constants, and it has a natural structure as
a Riemannian manifold. A curve in $\mathcal{K}$ can be represented by
a function $\psi(t,z)$ on $Z$, depending on an additional parameter
$t$. Here $t$ is naturally a real variable but following the usual
convention we let $t$ be complex, and assume $\psi$ independent of the
imaginary part  of $t$. Our domain $U$ is then a strip and metrics
$\phi$ satisfying $(i\ddbar\phi)^{n+1}=0$ correspond to geodesics in
$\mathcal{K}$.

 The basic strategy in \cite{Donaldson 1} and \cite{Donaldson 2} is
to approach the study of metrics on $\hat L$ - and hence Kähler metrics
on $Z$ with Kähler form in $c(\hat L)$- through the induced metrics
$Hilb({\phi})$ on 
$H^0(Z, \hat L)$, defined by
$$
\|u\|^2 =\int_Z |u|^2 e^{-p\phi}\omega^{\phi}_n,
$$
and similar metrics on $H^0(Z, \hat L^p)$ for $p$ large. 
Here $\omega_n=\omega^{\wedge n}/n!$ is the volume element on $Z$
induced by $\omega$. In our terminology,  \cite{Donaldson 1} and
\cite{Donaldson 2} thus deal with the bundle $F$.

\bigskip

In this paper we will follow a  different approach, using the
spaces $\hat E= H^0(Z,L\gr K_z)$ and the induced norms $H_\phi$
described above instead of $H^0(Z, L)$ and $Hilb({\phi})$. 

In \cite{Donaldson 2}, Donaldson introduced  certain, intimately
related,  functionals on
$\mathcal{K}$, $\mathcal{L}_p$ and $\tilde{\mathcal{L}}_p$. Stationary
points of $\tilde{\mathcal{L}_p}$ 
correspond to so called balanced metrics on $\hat L^p$. As $p$ tends
to infinity, $\tilde{\mathcal{L}}_p$ goes to the so called Mabuchi
functional, whose 
stationary points are precisely the metrics of constant scalar
curvature on $Z$. Donaldson's $\mathcal{L}$ functional can be defined
in terms of the metric on det($F$) induced by the metric on $F$. The
negativity of $F$, and hence det$(F)$, then give convexity of
$\mathcal{L}$ along geodesics (which is mentioned without proof in
\cite{Donaldson 2}). 

Here we  introduce the analogs of Donaldson's functionals
in our setting and show that they have largely the same properties,
with one important difference: In our setting, the $\mathcal{L}$
functional is defined with the opposite sign compared to
\cite{Donaldson 1} and \cite{Donaldson 2}, reflecting the different
estimates for the curvatures of the bundles $F$ and $E$.

Nevertheless our $\mathcal{L}$-functional  still converges to the
Mabuchi functional as $p$ goes 
to infinity. This fact is actually somewhat easier to verify in our
setting. We also show that the second derivative of the
$\tilde{\mathcal{L}}_p$-functional along a curve can be computed from the trace
of the curvature 
of our associated vector bundles. From this it follows that the second
derivative of the Mabuchi functional along a geodesic in $\mathcal{K}$
equals the limit
$$
\lim_{p\rightarrow \infty} \tau(p).
$$
From this we see that our asymptotic formula for the trace $\tau(p)$
is in fact equivalent to a known formula for the second derivative of
the Mabuchi functional, see   \cite{Donaldson 3}.

As is well known, this formula gives a ``formal proof'' of the
uniqueness properties of metrics 
of constant scalar curvature: Using the characterization of when the
above limit vanishes, we see that two metrics of constant scalar
curvature can be connected through the flow of a holomorphic vector
field, provided that they can be connected with a smooth geodesic in
the Riemannian space $\mathcal{K}$. This certainly does not give the
results in \cite{Donaldson 1}, \cite{Bando-Mabuchi},\cite{Mabuchi 1}
and \cite{Chen-Tian}, as the
existence of smooth geodesics is an open and even dubious
issue, (see \cite{Chen} and \cite{Chen-Tian} for the best results in
this direction), but at least it indicates a seemingly interesting
relation between our results 
on flatness of direct image bundles and the results of
\cite{Bando-Mabuchi},  \cite{Mabuchi 1} 
and \cite{Chen-Tian}. 

It is proved in \cite{Chen} that two metrics in $\mathcal{K}$ can
always be connected with a geodesic of class $C^{1,1}$. Phong and
Sturm have showed that this geodesic can be obtained as an almost
uniform limit of certain ``finite dimensional geodesics'' obtained
from Bergman kernels in $H^0(Z, L^p)$. In the last section we 
show that in our setting we can improve on this result somewhat, and
prove that Bergman kernels for the spaces $H^0(Z,L^p\gr K_Z)$
approximate the geodesic uniformly at the rate at least $p/\log
p$. This section is independent of most of the rest of the paper,
using only Theorem 2.1. Therefore a reader mainly interested in the
convergence of geodesics can go directly to the last section after
reading Theorem 2.1.

Part of the results in this paper were announced in \cite{B2}. As is
probably clear from the text, I am a novice in 
the study of extremal Kähler metrics, and I apologize  for
any omission in accrediting results properly. The motivation for this
paper lies not so much in the particular results as in showing the
link between $\dbar$-theory and this beautiful area.  It
should also be stressed that the list of references in this paper is far
from complete. Finally I would like to thank Robert Berman, Sebastien
Boucksom and Yanir Rubinstein for helpful and stimulating discussions
on these matters.

\section{ Positivity of direct images}

We start by giving the precise form of the result from \cite{B1} in
this particular setting, a trivial fibration with compact fibers over
a one-dimensional base. We will assume all the time that the metric
$\phi$ restricts to a positive metric on each fiber $\{t\}\times
Z$. Most of the time, but not always, we will also assume that $\phi$
defines a semipositive metric on the total space $X$. With
$\phi=\phi^0 +\psi$ given we put 
$$
c(\phi) = \psi_{t \bar t} -|\dbar_z\psi_t|^2_{\omega^t}.
$$
Here $\psi_t$ means $\partial\psi/\partial t$ and the subscript $\omega^t$
indicates that we measure the $(0,1)$-form $\dbar_z\psi_t$ with the
metric
$$
\omega^t =i\ddbar_z\phi|_t.
$$
This function  plays a major role in our estimates, and also in the
theory of variations of Kähler metrics. A short
calculation, see \cite{Mabuchi}, \cite{Semmes} or \cite{B1} shows that
it is related to the complex Monge-Ampere operator through the formula
\begin{equation}
c(\phi)idt\wedge d\bar t\wedge \omega^t_n = (i\ddbar \phi)_{n+1}.
\end{equation}
It is also, as proved in the references above, equal to the geodesic
curvature of the path defined by $\phi$ in the Riemannian manifold
$\mathcal{K}$. 

To state our formula for the curvature of $E$ it is convenient to
introduce yet another piece of notation. Let $f$ be a $\dbar$-closed
form of bidegree $(n,1)$ with values in $\hat L$ on $Z$. Given a
positive metric $\phi$ on $\hat L$ we can, by the Hörmander $\dbar
$-estimate solve $\dbar v=f$ with the $L^2$-estimate
$$
\|v\|^2\leq \|f\|^2,
$$
where the norms are defined using the Kähler metric $i\ddbar\phi$ on
$Z$. Let $v$ be the $L^2$-minimal solution and put
$$
e(f)=\|f\|^2 -\|v\|^2.
$$
Thus $e$ is a quadratic form in $f$, which by the Hörmander estimate is
positively semidefinite. If we develop $f$
$$
f=\sum f_j,
$$
where the $f_j$ are eigenforms of the $\dbar$-Laplacian with
eigenvalues $\lambda_j$,  one can easily verify that
$$
e(f)=\sum (1-\frac{1}{\lambda_j})\|f_j\|^2,
$$ 
but we will not use this. We next restrict our quadratic form $e$  to
$\dbar$-closed forms $f$ of the form 
$$
f=\dbar\mu\wedge u
$$
where $\mu$ is a smooth function and $u$ is a holomorphic $L$ valued
$(n,0)$-form. We then put
$$
A(\mu,u):= e(\dbar\mu\wedge u).
$$
$A$ is a quadratic form in both $u$ and $\mu$ separately, and it will
play a main role in the sequel. We will give a conjectural geometric
description of $A$ as the Chern curvature form of a vector bundle over a
certain infinite dimensional complex manifold in section 5.
\begin{thm} Let $u$ be an element in $E_t$. Then
$$
\langle \Theta^E u,u\rangle =\int_{\pi^{-1}(t)} c(\phi)[u,u]e^{-\phi}
+ A(\psi_t,u).
$$
\end{thm}

This result is implicit in \cite{B1}, but we shall indicate a simple
proof in this special situation. Take $t=0$. We extend $u$ to a
holomorphic section to $E$ near the origin in such a way that $D'u =0$
at 0, where $D'$ is the Chern connection on $E$. Then 
$$
\langle \Theta^E u,u\rangle =-\frac{\partial^2}{\partial t\partial\bar
  t}|_{t=0}\|u(t)\|^2_t .
$$
It is easily verified, cf e g \cite{B1}, that
$$
D'u = \Pi_{\text{holo}}(u_t -\psi_t u),
$$
with $\Pi_{\text{holo}}$ being the projection on the subspace of holomorphic
sections. Hence,  $D'u=0$ means that $v=(u_t -\psi_t u)$ is orthogonal
to the space of holomorphic forms. Therefore, for $t=0$, $v$ is the
$L^2$-minimal 
solution to the equation
$$
\dbar_z v=-\dbar_z\psi_t\wedge u=-f.
$$
Since
$$
\frac{\partial}{\partial t}\|u(t)\|^2_t= \langle u_t-\psi_t u,u\rangle,
$$
we get
$$
\frac{\partial^2}{\partial t\partial\bar t}\|u(t)\|^2_t |_{t=0}=
 \frac{\partial}{\partial \bar t}\langle(u_t-\psi_t u),u\rangle=
-\int\psi_{t \bar t} [u,u]e^{-\phi} +\int[v,v]e^{-\phi}.
$$
By definition
$$
\int[v,v]e^{-\phi}=\|v\|^2=\|f\|^2 -e(f),
$$
so
$$
\langle \Theta^E u,u\rangle=
$$
$$
=\int_{\pi^{-1}(t)} (\psi_{t \bar t}[u,u]-|f|^2)e^{-\phi}+e(f)
=\int_{\pi^{-1}(t)} c(\phi)[u,u]e^{-\phi}
+e(f),
$$
and the proof is complete. 

\bigskip

Since by Hörmander's theorem $e(f)$ and hence $A$ are always nonnegative
it follows immediately from Theorem 2.1 that
\be
\langle \Theta^E u,u\rangle \geq\int_{\pi^{-1}(t)} c(\phi)[u,u]e^{-\phi},
\ee
with equality if and only if $A(\psi_t,u)=0$. We also see that in case $\phi$
is assumed to be semipositive, then $\Theta^E$ is also semipositive
and can have a null-vector $u$ only if $c(\phi)=0$, and
$e(\dbar_z\psi_t\wedge u)=0$. By (2.1) the first condition means that $\phi$
satisfies the homogenuous Monge-Ampere equation 
$$
(i\ddbar\phi)^{n+1}=0.
$$
To understand the meaning of the second condition, we need to analyze
the quadratic form $e(f)$ a bit further.

Recall that the minimal solution $v$ to the equation $\dbar v=f$ can
be written
$$
v=\dbar^*\alpha
$$
for some $(n,1)$-form $\alpha$. This is simply because on a compact
manifold the range of $\dbar$ and $\dbar^*$ are closed so any form
orthogonal to the null-space of $\dbar$ lies in the image of the
adjoint operator. Moreover, by taking $\alpha$ orthogonal to the
kernel of $\dbar^*$, we can assume that $\dbar\alpha=0$, and $\alpha$
is then uniquely determined (as follows from (2.4) below). We will use
the Kodaira-Nakano identity in the following form.
\begin{equation}
\|\alpha\|^2 +\|v\|^2 +\|\dbar*\alpha\|^2 =2\langle
f,\alpha\rangle,
\end{equation}
where all norms are taken with respect to the Kähler metric given by
the curvature form of the metric on $\hat L$.
This slightly nonstandard formula follows from the more standard
\be
\|\alpha\|^2 +\|\dbar*\alpha\|^2= \|\dbar\alpha\|^2 +\|\dbar^*\alpha
\|^2,
\ee
if we add $\|\dbar^*\alpha\|^2=\|v\|^2$ to both sides, and use 
$$
\|\dbar^*\alpha\|^2=\langle f,\alpha\rangle.
$$

Now,
$$
2\langle f,\alpha\rangle=2\Re \langle
f,\alpha\rangle=\|f\|^2+\|\alpha\|^2 - \|f-\alpha\|^2.
$$
Inserting this in (2.3) and simplifying we find that
\be
e(f)=\|f-\alpha\|^2 +\|\dbar *\alpha\|^2 .
\ee
Thus, $e(f)=0$ if and only if $f=\alpha$ and $\dbar *\alpha=0$, so if
$e(f)$ vanishes we
must have $\dbar *f=0$. Conversely, if  $\dbar *f=0$, then $\dbar
\dbar^* f=f$ since 
$$
\dbar\dbar^* f=\dbar(-e^{\phi}\partial e^{-\phi}* f)=\omega\wedge
*f=f.
$$
Hence $\alpha=f$ and we see that $e(f)=0$. We have therefore proved
the next proposition.
\begin{prop}
$$
e(f)=0,
$$
i e equality holds in the Hörmander estimate for the $L^2$-minimal
solution to 
$$\dbar v=f,
$$
if and only if the
$(n-1,0)$-form $*f$ is holomorphic.
\end{prop}

\bigskip

For later use we also give a variant of the argument above that
leads to a  more precise statement.

\begin{prop}
$e(f)$ is the minimum of the quadratic expression
$$
q(f,g):=\|f-g\|^2 +\|\dbar*g\|^2
$$
where $g$ ranges over all ($L$-valued) $\dbar$-closed
$(n,1)$-forms. Hence, if we take  $g=f$, it follows that
$$
e(f)\leq \|\dbar*f\|^2.
$$
\end{prop} 
\begin{proof}
We already know that $q(f,\alpha)=e(f)$ if $\alpha$ is chosen so that
$\dbar\dbar^*\alpha=f$ and $\dbar\alpha=0$, so we need only prove that
for any $g$ as above $q(f,g)\geq e(f)$. Write $g=\alpha+\gamma$. We
claim that
$$
\langle\alpha-f,\gamma\rangle +\langle \dbar*\alpha,\dbar
*\gamma\rangle=0.
$$
This follows from polarizing (2.4); if we recall that
$\dbar\alpha=\dbar\gamma=0$ we get
$$
\langle\alpha,\gamma\rangle +\langle \dbar*\alpha,\dbar
*\gamma\rangle=\langle\dbar^*\alpha,\dbar^*\gamma\rangle=\langle
f,\gamma\rangle.
$$
Hence
$$
q(f,\alpha+\gamma)=q(f,\alpha)+\|\gamma\|^2+\|\dbar*\gamma\|^2\geq
e(f).
$$
\end{proof}

Recall that we are interested in when $e(f)=0$ when
$f=\dbar_z\psi_t\wedge u$. Let $V$ be the {\it complex gradient} of
$\psi_t$ on a fiber $\{t\}\times Z$. This is a vector field of type
$(1,0)$ defined by 
$$
\delta_V\omega^t=\dbar\psi_t,
$$
where $\delta_V$ means contraction of a form with a vector field. We
claim that $*(\dbar_z\psi_t\wedge u)=-\delta_V u$. To
see this, note that
$$
\omega^t\wedge u=0
$$
for degree reasons, so
$$
0=\delta_V(\omega^t\wedge u)=\dbar_z\psi_t\wedge u
+\omega^t\wedge\delta_V u.
$$
Hence
$$
*\delta_V u= -\dbar_z\psi_t\wedge u,
$$
which proves the claim.

Since $u$ is a holomorphic form, $\delta_V u$ is holomorphic if and
only if $V$ is a holomorphic vector field. This is clear outside of
the zeros of $u$, and since $V$ is automatically smooth, it must hold
everywhere.  We therefore see that for any
choice of $u$ in $E_t$ and 
$f=\dbar_z\psi_t\wedge u$, we have $e(f)=A(\psi_t,u)=0$ if and only if $V$, the
complex gradient of $\psi_t$, is a holomorphic vector field. Note that,
somewhat surprisingly, 
this condition is independent of $u$. We have
therefore proved the next theorem.
\begin{thm}
Equality holds, for some $u$, in the inequality
$$
\langle \Theta^E u,u\rangle \geq\int_{\pi^{-1}(t)} c(\phi)[u,u]e^{-\phi},
$$
if and only if $V$, the complex gradient of $\psi_t$ on the fiber
$Z_t=\{t\}\times Z$, is a holomorphic vector field. If $\phi$ is
semipositive, then $\Theta^E$ has a null vector in $E_t$ if and only
if $c(\phi)=0$ and $V$ is holomorphic on $Z_t$. In this case
$\Theta^E$ vanishes on  all of $E_t$.
\end{thm}

We shall next see that if the conditions of the previous theorem -
vanishing of $c(\phi)$ and holomorphicity of $V$ on $Z_t$ - are
satisfied for $t$ in an entire neighbourhood of the origin, then the variation
of the metrics on $Z$ comes from the flow of the holomorphic field
$V$, the complex gradient of $\psi_t$.

Let us first recall a few basic facts about real and complex vector
fields. If $v$ is a real vector field on $X$ then $v$ generates a
flow, a one-parameter group  $F_\tau$, of diffeomorphisms of $X$, satisfying
$$
\frac{d}{d\tau}g(F_\tau(x))=v(g)(F_\tau(x))
$$
for any smooth function $g$. The Lie derivative of a form $\alpha$ 
on $X$ with respect to $v$ is defined by 
$$
\mathcal{L}_v\alpha =\frac{d}{d\tau}|_{\tau=0}F_\tau^*(\alpha).
$$
By a classical formula of Cartan
$$
\mathcal{L}_v\alpha =(d\delta_v+\delta_v d)\alpha.
$$
If $V$ is a complex vector field of type $(1,0)$ we similarily let
the {\it complex Lie derivative}  be
$$
\mathcal{L}_V^{\C}\alpha =(\partial\delta_V+\delta_V \partial)\alpha.
$$
When $V$ is a holomorphic field, and $v=\Re V$  one easily checks that 
$$
\mathcal{L}_v=\Re \mathcal{L}_V^\C,
$$
so information about the complex Lie derivative of $V$ enables us to draw
conclusions about the flow of $v=\Re V$. 

Recall that the complex gradient of $\psi_t$ is defined by
$$
\dbar\psi_t=\delta_V\omega^t.
$$
Applying $i\partial$ to both sides we get (since $\partial \omega=0$)
$$
\frac{\partial}{\partial  t}\omega=i
\partial\delta_V\omega=i\mathcal{L}_V^\C\omega
$$
so we have a formula for the complex Lie derivative of $\omega$. To
use this formula, we need however to know that $V$ is a holomorphic
field, not just on each slice $X_t$, but also that it depends
holomorphically on $t$. That is the object of the next lemma.

\begin{lma} Assume $V$, the complex gradient of $\psi_t$ is
  holomorphic on each slice $X_t=Z\times\{t\}$ for all $t$ in $U$. Let $V_{\bar
    t}=\partial V/\partial 
  \bar t$. Then $V_{\bar t}$ is the complex gradient of $c(\phi)$. In
  particular, if $c(\phi)$ is constant on slices then $V$ is a holomorphic field
  on all of $X=U\times Z$. 
\end{lma} 
\proof

Recall that
\be
\frac{\partial}{\partial  t}\omega=i
\partial\delta_V\omega=i\mathcal{L}_V^\C\omega,
\ee
where $\mathcal{L}_V^\C$ is the complex Lie derivative with respect to
$V$. Taking conjugates it follows that
$$
\frac{\partial}{\partial  \bar t}\omega=-i
\dbar\delta_{\bar V}\omega.
$$
Differentiating the defining relation
$$
\dbar\psi_t=\delta_V\omega^t
$$
with respect to $\bar t$ we get
$$
\dbar_z\psi_{t \bar t}=\delta_{V_{\bar t}}\omega
+\delta_V\frac{\partial}{\partial  \bar t}\omega =
$$
$$
\delta_{V_{\bar t}}\omega-i\delta_V\dbar\delta_{\bar V}\omega=
\delta_{V_{\bar t}}\omega+i\dbar\delta_V\delta_{\bar V}\omega=
\delta_{V_{\bar t}}\omega+\dbar |V|^2.
$$
Thus
$$
\delta_{V_{\bar t}}\omega=\dbar(\psi_{t \bar t}-|V|^2),
$$
and the lemma follows since the norm of $V$ is equal to the norm of
$\dbar\psi_t$. 
\qed
\begin{thm}Assume the curvature $\Theta^E$ is degenerate on $E_t$ for
  all $t$ in $U$. Then
$$
\omega^t= S_t^*(\omega^0),
$$
where $S_t$ is the flow of some holomorphic vector field on $Z$. It
follows that
$$
\phi(t,z)=\phi(0,S_t(z)) +\psi(t).
$$
\end{thm}
\proof
By Theorem 2.3 and Lemma 2.4 the hypothesis implies that $V$, the
complex gradient of $\psi_t$ is holomorphic with respect to both $t$
and $z$. Let 
$$
T_s(t,z)
$$
be the flow of $iV$; note that it acts fiberwise on $X$. Let
$$
U_s(t,z)
$$
be the flow of $W:=iV-\partial/\partial t$.
Then
$$
U_s(t,z)=(t-s,T_s(t,z)).
$$
By (2.7), the (complex) Lie derivative of $\omega$ with respect to $W$
vanishes, so $ U_s^*(\omega)=\omega$. Taking $s=t$ it follows that 
$\omega^0=T_t^*(\omega^t)$ and the theorem follows with $S$ equal to the
inverse of $T$. (The last part follows since
$\phi(t,z)-\phi(0,S_t(z))$ is pluriharmonic, hence constant, on fibers.) 

\qed

\section{Negativity of direct image bundles.}

In this section we give, for comparison, a theorem that very roughly
corresponds to Theorem 2.1, for the (trivial) vector bundle $F$ with fiber
$$
\hat F=F_t= H^0(Z, L),
$$
instead of $\hat E=H^0(Z,K_Z\gr L)$. The difference is thus that we do
not take tensor products with the cannonical bundle here.  We
then  use the metric
$$
\|u\|^2_t:=\int_{\pi^{-1}(t)}|u|^2 \ephi\omega_n,
$$
where $\omega=(i\ddbar)_z\phi$, on $F$. We shall see that
we then instead get an estimate from above 
of the curvature. 

\begin{thm}Let $\phi$ be a smooth metric on $L$ over $Z\times U$
  satisfying
$$
\omega^\phi:=i\ddbar\phi\geq 0.
$$
Let $u$ be an element in $F_t$. Then the curvature $\Theta^F$ of $F$
equipped with the metric described above satisfies
$$
\langle\Theta^Fu, u\rangle \leq (n+1)\int_{\pi^{-1}(t)}|u|^2
c(\phi)\ephi\omega_n.
$$
\end{thm}
\begin{proof} Let $u$ be a holomorphic section to $F$. We start by
  computing
$$
i\ddbar \|u\|^2_t.
$$
For this we note that $\|u\|^2$ is the push-forward of the form 
$$
R=|u|^2\ephi\omega^\phi_n
$$
under $\pi$. Hence, if we denote by $\partial^\phi$ the twisted
derivative $e^\phi\partial\ephi$,
$$
i\ddbar \|u\|^2_t=\pi_*(i\ddbar R)=
$$
\be
\pi_*(i\partial^\phi
u\wedge\overline{\partial^\phi u}\wedge\omega^\phi
\ephi)-\pi_*(|u|^2\omega^\phi\wedge \omega^\phi_n\ephi).
\ee
Here we have used the commutator rule
$$
\dbar\partial^\phi +\partial^\phi\dbar =\ddbar\phi.
$$
The first term in the right hand side of (3.1) is nonnegative, so using
$$
\omega^\phi\wedge \omega^\phi_n=(n+1)c(\phi)idt\wedge d\bar t\wedge \omega^\phi_n
$$
we get that
\be
i\ddbar \|u\|^2_t\geq -(n+1)\int_{\pi^{-1}(t)}|u|^2
c(\phi)\ephi\omega_n idt\wedge d\bar t. 
\ee
On the other hand
$$
i\ddbar \|u\|^2_t=\langle D'u,D'u\rangle -\langle \Theta^F u,
u\rangle,
$$
if $D'$ is the $(1,0)$-part of the connection on $F$.
If we combine this formula with (3.2) we see that
$$
\langle \Theta^F u,u\rangle\leq \langle D'u,D'u\rangle +
(n+1)\int_{\pi^{-1}(t)}|u|^2 c(\phi)\ephi\omega_n idt\wedge d\bar t.
$$
Since we can make the first term in the right hand side vanish at any
given point by choosing the section $u$ so that $D'u=0$ at that given
point, it follows that
$$
\langle\Theta^Fu, u\rangle \leq (n+1)\int_{\pi^{-1}(t)}|u|^2
c(\phi)\ephi\omega_n.
$$
\end{proof}
\begin{cor}Assume, in addition to the hypotheses in the previous
  theorem that $\omega^\phi_{n+1}=0$, so that $c(\phi)=0$. Then
  $\Theta^F\leq 0$
\end{cor}

\bigskip

\noindent{\bf Remark:} It might seem surprising (and even confusing)
that we get different signs for the curvature of the very closely
related bundles $E$ and $F$. This can be clarified somewhat if we
consider the real variable analogs of Theorems 2.1 and 3.1. In this
analogy a metric on $L$ over $Z\times U$ with nonnegative curvature
corresponds to a convex  function on $\R^n\times U$, with $U$ now an
interval in $\R$.  As explained in \cite{B1}, Theorem 2.1 
corresponds to {\it Prekopa's Theorem}, which states that the function

$$
-\log\int_{\R^n}e^{-\phi(x,t)}dx
$$
is convex. On the other hand, Theorem 3.1 corresponds to the fact that

$$
\log\left(\int e^{-\phi(x,t)}\text{det}(\phi_{x_j x_k})dx\right)
$$
is convex, if the determinant of the full hessian of $\phi$ vanishes.

\section{Asymptotics .}

As in the introduction we next let, for $p$ a positive integer, $E(p)$ be the
vector bundle defined in 
the same way as $E$, but replacing $L$ by $L^p$. It follows
immediately from Theorem 2.3 that if the curvature $\Theta^p$ of
$E(p)$ is degenerate at $t$ for $p$ equal to some $p_0$, then actually
$\Theta^p$ 
vanishes completely at $t$ for any $p=1,2....$. We shall now see that
we can draw a similar conclusion if $\Theta^p$ vanishes in the limit
as $p$ tends to infinity. This requires an asymptotic study of the two
terms in our curvature formula from Theorem 2.1.

\bigskip

\noindent  The main point in our asymptotic study of $\Theta^p$ is the
next asymptotic formula for the trace of the quadratic form
$A(\mu,\cdot)$ introduced in section 2. 
\begin{thm} Let $A_p$ be the quadratic form $A$ with $L$ replaced by
  $L^p$. Denote by $d_p$ the dimension of $H^0(Z,K_Z\otimes L^p)$ and
  by Vol$(Z)$ the volume of $Z$ with respect to the metric $\omega$. Then
$$
\lim_{p\rightarrow \infty} \frac{1}{d_p}\mathbf{tr}
A_p(p\mu,\cdot)=\frac{1}{\text{Vol}(Z)}\int_Z |\dbar 
V_\mu|^2 \omega_n.
$$
\end{thm}
In the proof we will use the asymptotic expanison for Bergman kernels
of Tian-Catlin-Zelditch that we recall in the next subsection.

\subsection{The Tian-Catlin-Zelditch asymptotic formula for Bergman
  kernels.}

Fix $p$ for the moment, and let $u_j$ be an
orthonormal basis for $E(p)=$. Then
$$
\sum [u_j,u_j]e^{-p\phi}=:B_pe^{-p\phi}
$$
is the Bergman  form for $E(p) =H^0(Z,K_Z\otimes L^p)$.
By  e g \cite{Zelditch},
 there are constants $a_j$ so that
\be
B_{p\phi}e^{-p\phi} =p^n\left(a_0 +\frac{a_1}{p}S
  +O(\frac{1}{p^2}) \right)\omega^\phi_n.
\ee
Here $a_0$ and $a_1$ are positive constants, and by the computation of
the second term due to Lu, \cite{Lu}, $S=S_{\omega^{\phi}}$ is the scalar
  curvature of the metric with Kahler form $\omega^\phi$. 
 Normally the expansion
is given for the Bergman kernel for the space of global sections to
$L^p$, but the analogous formula for our Bergman kernel follows from
this. More precisely, it follows from the known expansions for $L^p\gr
F$, where $F$ is another line bundle with fixed metric - our case
corresponding to $F=K_z$ with metric induced by $\omega^\phi$, see e g
\cite{Berman-Berndtsson-Sjöstrand}.

\bigskip

\noindent Since 
$$
\frac{ d_p}{p^n} =b +c/p +O(\frac{1}{p^2})
$$
where $b$ is again a positive constant it follows from 2.7 that
\be
\frac{B_{p\phi}e^{-p\phi}}{d_p} =\left(a'_0
  +\frac{a'_1}{p}(S-\hat S)
  +O(\frac{1}{p^2}) \right)\omega^\phi_n
\ee
for certain constants $a_0'$, $a_1'$ and $\hat S$. 
The left hand side here integrates to 1, so it follows that
$a'_0=1/\text{Vol}(Z)$. Moreover, the second term in the right hand
side must integrate to 0, so $\hat S$ is the average of $S$ over $Z$.

\subsection{ Proof of Theorem 2.7}
We shall first prove that
$$
\limsup_{p\rightarrow \infty}\frac{1}{d_p} \mathbf{tr}
A_p(p\mu,\cdot)\leq\frac{1}{\text{Vol}(Z)}\int_Z |\dbar 
V_\mu|^2 \omega_n.
$$
This is a relatively simple consequence of Proposition 2.3, together 
with the formula for the asymptotic expansion of Bergman forms.  
Recall that, by Proposition 2.3, 
$$
A_p(p\mu,u)\leq \|p\dbar*_p\dbar\mu\wedge u\|_p^2
=\|\dbar*\dbar\mu\wedge u\|^2   .
$$ 
Here the $p$ in the subscript indicate norms and *-operators with
respect to the metric $p\omega$, and norms and *-operators without
subscripts are taken with respect to $\omega$. (Notice that $p*_p=*$
and that the norm of a form of total degree $n$ does not change when
we multiply the metric with $p$. ) But, it is easily verified that
$$
\|\dbar*\dbar\mu\wedge u\|^2 =\int |\dbar V_\mu|^2[u,u].
$$
This means that the trace of $A_p$ is dominated by
$$
\int |\dbar V_\mu|^2 B_{p\phi},
$$
which after divison with $d_p$ tends to
$$
\frac{1}{\text{Vol}(Z)}\int |\dbar V_\mu|^2 \omega_n
$$
by the asymptotic formula for $B_{p\phi}$. This completes the proof of
the upper estimate. We next turn to the estimate from below of the
trace, i e we estimate $\dbar V_\mu$ from above.

\bigskip

\noindent For this we shall estimate
$$
I=\int (\dbar V_\mu, \xi)\omega_n\text{Vol}(Z)
$$
where $\xi$ is a $(0,1)$ with vector field coefficients. By the
asymptotic expansion for Bergman forms this is the limit of 
$$
I_p=\int  (\dbar V_\mu, \xi)\sum [u_j,u_j]/d_p
$$
as $p$ tends to infinity, if $(u_j)$ is an orthonormal basis for
$H^0(Z, K_Z\gr L^p)$. 
But
$$
I_p=\sum \langle\dbar\delta_{V_\mu}u_j,\delta_\xi u_j\rangle/d_p=
\sum \langle\delta_{V_\mu}u_j,\dbar^*\delta_\xi u_j\rangle/d_p.
$$
Recall that
$$
A_p(p\mu,u_j)=\|p\dbar\mu\wedge u_j-p\alpha_j\|_p^2 +\|\dbar*_p\alpha_j\|_p^2
$$
where $\alpha_j$ are $(n,1)$-forms. Translating to norms with respect
to $\omega$ this equals
$$
p\|\dbar\mu\wedge u_j-\alpha_j\|^2 +\|\dbar*\alpha_j\|^2.
$$
In particular
$$
\sum \|\dbar\mu\wedge u_j-\alpha_j\|^2/d_p=\sum \|\delta_{V\mu}
u_j-*\alpha_j\|^2/d_p \leq \frac{1}{p}\mathbf{tr}A_p(p\mu,\cdot)/d_p.
$$
Hence, up to an error of size $1/p$, $I_p$ equals
$$
\sum \langle*\alpha_j,\dbar^*\delta_\xi u_j\rangle/d_p=
\sum \langle\dbar*\alpha_j,\delta_\xi u_j\rangle/d_p,
$$
which by Cauchy's inequality is dominated by
$$
(\sum \|\dbar*\alpha_j\|^2/d_p)^{1/2}(\sum \|\delta_\xi
u_j\|^2/d_p)^{1/2}\leq
$$
$$
\leq
(\mathbf{tr}A_p(p\mu,\cdot)/d_p)^{1/2}( \int |\xi|^2 B_{\phi}/d_p)^{1/2}.
$$
Therefore
$$
\frac{1}{\text{Vol}(Z)}\int (\dbar V_\mu, \xi)\omega_n\leq
\liminf(\mathbf{tr}A_p(p\mu,\cdot)/d_p)^{1/2}( \|\xi\|^2/\text{Vol}(Z))^{1/2},
$$
so
$$
\frac{1}{\text{Vol}(Z)}\int |\dbar
V_\mu|^2\omega_n\leq\liminf \frac{1}{d_p}\mathbf{tr}A_p(p\mu,\cdot).
$$
This completes the proof of Theorem 2.7. \qed 

\subsection{Asymptotic behaviour of the curvature}

In the previous subsection we have used the asymptotics of Bergman
forms to describe the asymptotic behaviour of the trace of $A_p$. The
behaviour of the first term in the formula for the curvature $\Theta^p$ 
as $p$ goes to infinity also follows from the Tian-Catlin-Zelditch
formula. This term is
$$
\int c(p\phi)[u,u],
$$
so its trace is
$$
\int pc(\phi)B_{p\phi}.
$$
The results recalled in subsection 2.1 therefore lead to the next
theorem.
\begin{thm} As $p$ goes to infinity
$$
\frac{1}{d_p}\mathbf{tr} \Theta^p=
$$
$$
= (p\int c(\phi)\omega^{\phi}_n
+\int c(\phi)(S-\hat S)\omega^\phi_n +\int |\dbar
V_{\psi_t}|^2\omega^\phi_n)/\text{Vol}(Z) + o(1).
$$
In particular
$$
\liminf \frac{1}{d_p}\mathbf{tr} \Theta^p=0
$$
(if and) only if $c(\phi)=0$ and $V_{\psi_t}$ is a holomorphic vector
field.
\end{thm}

\section{Convexity on the space of Kähler metrics}

We now return to the space $\mathcal{K} $ of Kähler potentials on
$Z$. Fixing a base metric $\phi_0$, any other metric $\phi$ can be
written
$$
\phi=\phi_0 +\psi
$$
where $\psi$ is a function on $Z$. The space $\mathcal K$ is therefore
affine and its tangent space is the space of (say smooth) real valued
functions on 
$Z$. Following \cite{Mabuchi} and \cite{Semmes} one defines a metric
on the tangent space at $\phi$ by
$$
\|\psi\|^2 =\int |\psi|^2 \omega^\phi_n,
$$
where $\omega^\phi=i\ddbar\phi$. This gives $\mathcal K$ the structure
of a Riemannian manifold and it is proved in the references above that
the geodesic curvature of a path $\phi(t,z)=\phi_0(z)+\psi(t,z)$ is
given by
$$
c(\phi) = \psi_{t \bar t}- |\dbar_z\psi_t|^2_{\omega^\phi}
$$
Here we follow the same convention as in the previous paragraph that
we let $t$ be complex and let $\psi$ be independent of the argument of
$t$. By formula (2.1) $\phi$ defines a geodesic in $\mathcal{K}$
precisely when $\phi$ satisfies the homogenuous Monge-Ampere equation
with respect to both variables $t$ and $z$.

Let us now consider a (smooth) function, $\mathcal F$ on $\mathcal K$. 
Then 
$$
\frac{\partial}{\partial t}\F(\phi_t)=\F'.\phi_t.
$$
To write second order derivatives we have two possibilities. The
simplest alternative would be to write
$$
\frac{\partial^2}{\partial t\partial \bar t}\F(\phi_t)=\F'.\phi_{t \bar
  t} + D^2\F(\phi_t,\phi_{\bar t}).
$$
Then $D^2\F$ is the Hessian of $\F$ with respect to the affine
structure on $\K$. The second possibility is to write
\be
\frac{\partial^2}{\partial t\partial \bar t}\F(\phi_t)=\F'.c(\phi)
+\F''(\phi_t,\phi_{\bar t}).
\ee
and let this formula  define the Hessian of $\F$ as
a quadratic form on the tangent space to $\K$. This is the Hessian of
$\F$ determined by the Riemannian structure on $\K$, and it is this
form of the Hessian that we will use. Then $\F$ is said to be
convex on $\K$ if $\F''$ is positively semidefinite, so that in
particular the restriction of $\F$ to any smooth geodesic
is convex. Notice however that the convexity is defined independently
of the existence of smooth geodesics.

One classical example of such a function is $I$, defined by
$$
I'.\phi_t=\int\phi_t\omega^\phi_n.
$$
It is well known, and not hard to compute, that
$$
\frac{\partial}{\partial \bar t}I'.\phi_t=\int c(\phi)
\omega^\phi_n=I'.c(\phi),
$$
so $I''=0$ and  $I$ is linear along geodesics. 

In \cite{Donaldson 2} Donaldson introduced another functional, $\mathcal{L}$,
in the following manner. Let $H$ be the space $H^0(Z,\hat L)$ of global
holomorphic sections to $\hat L$ over $Z$. Let $(h_j)$ be some fixed basis for
$H$ and let for any $\phi$ in $\K$, $\mathcal{L}(\phi)$ be the logarithm of the
determinant of the matrix of the metric $Hilb(\phi)$ expressed in the
given basis. Here the metric $Hilb(\phi)$ is defined by
$$
\|h\|^2 =\int |h|^2e^{-\phi}\omega^\phi_n.
$$
Note that $\mathcal{L}$ depends on the basis chosen, but only up to an additive
constant. Equivalently, Donaldson's $\mathcal{L}$-functional can be
defined as the logarithm of the square of the norm of some fixed
constant section of 
the line bundle det($F$), where $F$ is the vector bundle discussed in
section 3 (the notion of a constant 
section makes sense since $F$ and hence det($F$) are trivial
bundles). Corollary 3.2 therefore implies that
$\mathcal{L}(\phi(\cdot,t))$ is a subharmonic function of $t$ if
$\phi(\cdot,t)$ is any curve in $\mathcal{K}$ satisfying
$c(\phi)=0$. Hence, in particular, $\mathcal{L}$ is convex along
geodesics.Notice however that $\mathcal{L}$ has no apparent convexity
property along curves satisfying $i\ddbar\phi\geq 0$ or  $i\ddbar\phi\leq 0$. 

In this paper we will change the setup and notation  and let 
$\mathcal{L}$ be the 
similarily defined  functional, but {\it with opposite sign}, and
using the space $\hat 
E= H^0(Z,L\gr K_Z)$ instead 
of $H$, and  our metrics $H_\phi$ instead of $Hilb(\phi)$. Hence, in
our setup, $\mathcal{L}$ is the {\it negative} of the logarithm of the
squared norm of 
some fixed constant section of the line bundle det($E$), discussed in
section 2. Here the norm on det($E$) is of course the norm induced by
our norms $H_\phi$ on $E$.  

Thus $\mathcal{L}'$ is the negative of the 
connection
form on det($ E$) (with respect to the constant section chosen). Let
$(u_j)$ be an orthonormal basis of 
$E$ for one fixed $t$. Then  the connection form on det($E$)
is equal to  the trace of the 
connection, $\theta$  on $E$.
$$
\mathbf{tr}\theta=\sum\langle \theta u_j,u_j\rangle=-\sum\int \phi_t
[u_j,u_j]e^{-\phi}=-\int \phi_t B_\phi e^{-\phi},
$$ 
where $B_\phi$ is the Bergman kernel. We have thus proved the first
part of the following lemma.

\begin{lma}
$$
\mathcal{L}'.\phi_t=\int \phi_t B_\phi e^{-\phi}
$$
and
$$
\frac{\partial^2}{\partial t\partial \bar t}\L(\phi_t)=\mathbf{tr} \Theta^E.
$$
Hence
$$
\mathcal{L}''(\mu,\mu)=\mathbf{tr} A_{\mu}(\cdot,\cdot).
$$
\end{lma}
The second part of the lemma follows since the Laplacian of the
logarithm of the square norm of a holomorphic (in this case constant) section
to $\det E$ is the negative of the curvature of $\det E$, which is
equal to the trace of the curvature of $E$. The last part then follows
from the definition of the Hessian, formula (3.1), and Theorem 2.1.

Note that, by Lemma 3.1,  $\mathcal{L}''$ is nonnegative along
geodesics and vanishes 
only if $V$, the complex gradient of $\psi_t=\phi_t$ is a holomorphic
vector field on $Z$. Hence $\mathcal{L}$ is convex along geodesics, and also
convex along any curve such that $c(\phi)\geq 0$. 

Following \cite{Donaldson 2} we next let
$$
\tilde{\mathcal{L}}=\frac{1}{d}\mathcal{L}- \frac{1}{\text{Vol}}I,
$$ 
where $d$ is the dimension of $\hat E$ and $\text{Vol}$ is the volume
of $Z$. Then $\tilde{\mathcal{L}}'.\phi_t=0$ if $\phi_t$ is constant on $Z$ for
some $t$ and $\tilde{\mathcal{L}}$ is also convex along geodesics (since $I$ is
linear). 
The stationary points of $\tilde{\mathcal{L}}$ are in this setting the points
$\phi$ in $\K$ such that 
$$
B_\phi e^{-\phi} =\frac{d}{\text{Vol}}\omega^\phi_n.
$$
This follows immediately from 
$$
\tilde{\mathcal{L}}'.\mu=\int\mu(\frac{B_\phi
  e^{-\phi}}{d}-\frac{\omega^\phi_n}{\text{Vol}}).
$$
We will refer to these points as balanced, but note that this term now
has a slightly different meaning from what it has in \cite{Donaldson 2}.

Finally we give also the definition of the {\it Mabuchi functional},
$\mathcal{M}$. It is determined (up to a constant) by the formula for its
derivative
$$
\mathcal{M}'.\mu=\int\mu(S_{\omega^\phi} -\hat S_{\omega^\phi})\omega^\phi_n,
$$
where $S_{\omega^\phi}$ is the scalar curvature of the metric
$\omega^\phi$, and $\hat S_{\omega^\phi}$ is its average. 
Its critical points are precisely the metrics $\phi$ in $\K$ that
induce metrics of constant scalar curvature on $Z$.
\begin{prop} Let $\phi_0$ and $\phi_1$ be balanced points. Assume they
  can be joined by a smooth geodesic. Then 
$$
\omega^{\phi_0}= S^*(\omega^{\phi_1})
$$
where $S$ is the time 1 map of some holomorphic vector field on $Z$.
\end{prop}
To see this, let $\phi$ be the geodesic and consider the restriction
of $\tilde{\mathcal{L}}$ to the geodesic. Since the end points are balanced, the
derivative of $\tilde{\mathcal{L}}$ vanishes at the end points, and since
$\tilde{\mathcal{L}}$ is convex it must be constant. Since $I$ is
linear $\mathcal{L}''$ 
also vanishes so
$\Theta^E$ is zero by Lemma 5.1. The proposition then follows from
Theorem 2.4. 

We will now apply the same reasoning to $L^p$, where $p$ tends to
infinity. Define $\mathcal{L}_p$ the same way as $\mathcal{L}$, but
replacing $L$ by 
$L^p$. Put
\be
\tilde{\mathcal{L}}_p=\frac{1}{d(p)}\mathcal{L}- \frac{1}{\text{Vol}}I.
\ee
Then
$$
\tilde{\mathcal{L}}'.\mu=\int \mu \sigma_p
$$
where
$$
\sigma_p=\frac{B_{p\phi} e^{-p\phi}}{d(p)}-\frac{\omega^\phi_n}{\text{Vol}}.
$$
It follows (cf formula 4.2) that, as $p$ tends to infinity, $p\sigma_p$
tends to (a constant times) 
$$
(S_\phi-\hat S_\phi)\omega^\phi_n
$$
where $S_\phi$ is the scalar curvature of the metric $\omega^\phi$ and
$\hat S_\phi$ is the average of $S_\phi$. 
By definition of the 
Mabuchi functional, $\mathcal{M}$, this means that
$\tilde{\mathcal{L}}_p'$ tends to  
$\mathcal{M}'$. The next result says the second derivatives also
converge. The formula for the second derivative of the Mabuchi
functional in this proposition can be found in \cite{Donaldson 3}.

\begin{prop}

$$
\lim_{p\rightarrow\infty}\tilde{\mathcal{L}_p}''(\mu,\mu)=
\frac{1}{\mathbf{Vol}(Z)}  \int|V_\mu |^2\omega_n =\mathcal{M}''(\mu,\mu).
$$

\end{prop}
\proof
The first equality  follows directly from Lemma 5.1 and Theorem 4.1, and
the remaining equality follows of course since we know that the
first derivatives of $\tilde{\mathcal{L}_p}$ converge to  $\mathcal{M}'$.
\qed

{\bf Remark:} Conversely, Theorem 4.1 follows from the formula for
second derivative of the Mabuchi functional, if we use that $A$ is the
Hessian of $\tilde{\mathcal{L}}$ and that $\tilde{\mathcal{L}}_p$
converges to the Mabuchi functional. We have included the proof in
section 4, since it gives an approach to the Hessian of the Mabuchi
functional, using only data that can a priori be of very low
regularity.

In particular, the Mabuchi functional is convex along geodesics, and
even strictly convex if $Z$ has no nontrivial holomorphic vector
fields.
Recall that the critical points of $\mathcal{M}$ are precisely the Kähler
potentials such 
that $\omega^\phi$ has constant scalar curvature.
The same argument that lead to Proposition 3.2 now gives the next
theorem.
\begin{thm} Let $\phi_0$ and $\phi_1$ be  points in $\K$ that define
  metrics of constant scalar curvature on $Z$. Assume they
  can be joined by a smooth geodesic. Then 
$$
\omega^{\phi_0}= S^*(\omega^{\phi_1})
$$
where $S$ is the time 1 map of some holomorphic vector field on $Z$.
\end{thm}

\subsection{ An interpretation of the form $A$ (?)}

Our  form 
$$
A(\mu,u)=e(\dbar\mu\wedge u)
$$
 is defined on the tangent
bundle of $\mathcal{K}$ times the complex vector space 
$$
\hat E= H^0(Z,K_Z\gr L)
$$
and it is quadratic in both arguments so it resembles a curvature form
on a vector bundle 
over $\mathcal{K}$ with fiber $\hat E$. If we try to make this idea
more precise we must first find a connection. Clearly, the trivial
bundle over $\mathcal{K}$ with fiber $\hat E$ has a natural metric, if
we define the norm over a point $\phi$ in $\mathcal{K}$ to be just
$H_\phi$. A metric does not in itself induce a connection however. It
is therefore natural to consider complexifications of
$\mathcal{K}$, since holomorphic bundles over complex manifolds have
canonical connections . 

Let us therefore {\it assume} that we have a complex manifold
$\tilde{\mathcal{K}}$, containing $\mathcal{K}$ as a totally real
submanifold, together with a projection map
$$
\pi:\tilde{\mathcal{K}}\rightarrow \mathcal{K},
$$
having the property that any geodesic in $\mathcal{K}$ lifts to a
holomorphic curve in $\tilde{\mathcal{K}}$. Define a trivial
vector bundle $\mathcal{E}$ over $\tilde{\mathcal{K}}$ with fiber $\hat
E$ and the 'tautological' Hermitian metric $H_{\pi(z)}$ on
$\mathcal{E}_z$. We claim that $A(\mu,u)$ must then be the Chern
curvature form
$$
\langle\Theta^{\mathcal{E}}_{(\mu,i\mu)}u,u\rangle.
$$
To see this, take a piece of a geodesic curve through a point $\phi_0$
in $\mathcal{K}$ with tangent vector $\mu$ at $\phi_0$, and lift it to
a complex curve through a point $z_0$ in $ \tilde{\mathcal{K}}$. Then
$\mathcal{E}$ restricted to this complex curve is a vector bundle  $E$
of the type considered in section 2. By Theorem 2.1
its curvature is given by
$$
\langle\Theta^E u,u\rangle= A(\mu,u)
$$
(the first term in the curvature formula disappears for a
geodesic). But the Chern curvature of the bundle restricted to a
complex curve is the restriction of the Chern curvature of the full
bundle, proving our claim.

This interpretation hinges of course on the existence of a
complexification of $\mathcal{K}$ with the properties above. As
pointed out to me by Yanir Rubinstein, the at least 'moral' existence
of such a complexification is an important motivation for the
constructions in \cite{Semmes} and \cite{Donaldson 3}. One can also
compare to the work of Lempert and Szoke, \cite{L-Sz}, who show the
existence of a complex structure with the properties we require on a
neighbourhood of the zero section in the tangent bundle of any {\it
  finite dimensional} Riemannian manifold. At any rate, if
$\tilde{\mathcal{K}}$ exists we also see that
Hörmander's $L^2$-estimates imply (and are almost equivalent to) that
$\mathcal{E}$ has nonegative curvature, and that the curvature is
strictly positive if $Z$ has no nontrivial holomorphic vector fields.

\subsection{'Finite dimensional geodesics'.}

A Hermitian metric $H$ on the trivial vector bundle $E=U\times \hat E$ over $U$
is a complex curve in the space of Hilbert norms on $\hat E$.
The latter is a symmetric space (cf \cite{Semmes} and
\cite{Donaldson 1}) and it turns out that geodesics in this space
correspond exactly to metrics with zero curvature (at least if the
metric depends only on the real part of $t$). Let us call a metric
with semipositive curvature a {\it subgeodesic}, and in the same way
call a complex curve in $\mathcal{K}$ a subgeodesic if
$i\ddbar\phi\geq0 $. With this terminology, Theorem 2.6 says that if
$\phi$ is a subgeodesic in $\mathcal{K}$, then the corresponding
metric $H_\phi$ on $E$ is also a subgeodesic. Moreover, $H_\phi$ is a
geodesic, i e the curvature vanishes, if and only if $\phi$ arises
from one fixed metric on $\hat L$ via the flow of a holomorphic vector
field. 
 
Conversely, any metric on a vector bundle $F=U\times H^0(Z,L)$ induces a
corresponding curve in the space of metrics on $L$ by taking the
Bergman kernels of the norms on any fiber (this map from metrics on
$U\times H^0(Z,L)$ to $\mathcal{K}$ is called in the Fubini-Study
map, see \cite{Donaldson 2}). Assuming that $\hat L$ is sufficiently
positive so that the Bergman metric is nondegenerate on fibers, we 
have here a completely symmetric situation: Subgeodesics on $F$ map to
subgeodesics in $\mathcal{K}$, and the image is a geodesic if and only
if the metric on  on $F$ arises from one fix Hilbert norm on $\hat F$
via the flow of a holomorphic vector field. The first part of this
claim is clear from the explicit form of the Bergman kernels (see
beginning of next section or \cite{Phong-Sturm}). The second part
follows from an analysis of the foliation by complex curves that the
homogenuous Monge-Ampere equation 
$$
(i\ddbar\log B)^{n+1}=0
$$
induces, but we omit the proof.

\section{Approximation of geodesics}

In the previous section we did not use the full curvature estimate in
Theorem 2.1, but only the estimate of the trace of the curvature that
follows from it. In this section we shall use the full curvature
estimate to show that a recent result of Phong and Sturm,
\cite{Phong-Sturm} on approximation of geodesics can be sharpened a
bit in our setting.

Recall that any $\phi$ in $\K$ induces a Hilbert norm, $H_\phi$, on 
$\hat E$. Let $M$ be the space of all Hilbert norms on $\hat E$ and
let $H_t$ be a curve in $M$, where $t$ in $U$ is a complex parameter. Then
$H_t$ defines an Hermitian structure on our vector bundle $E$ over $U$.
We say that $H_t$ is flat if this Hermitian structure is flat, i e if
it has vanishing curvature. In case $H_t$ is independent of the
argument of $t$, this means precisely that the corresponding real
curve is a geodesic in the symmetric space of all Hermitian norms on
$\hat E$, see \cite{Donaldson 1} and \cite{Phong-Sturm}. 

As in \cite{Phong-Sturm}, we note that any two points, $H_0$ and $H_1$
in $M$ can be joined by a flat curve: Choose an orthonormal basis of
$H_0$, $(u_j)$ that diagonalizes $H_1$, so that
$$
\langle u_j,u_k\rangle_{H_1}=\delta_{j k} e^{2\lambda_j}.
$$
Then $H_t$, defined by
$$
\langle u_j,u_k\rangle_{H_t}=\delta_{j k} |t|^{2\lambda_j},
$$
where $\log|t|$ ranges from 0 to 1, is a flat curve, joining $H_0$ and $H_1$.

For any curve in $M$ we get a curve of metrics on $L\gr K_z$ by taking
the logarithm of $B_t$, the Bergman kernel for $\hat E$ with the
metric $H_t$. Here $B_t$ is defined by
$$
B_t=\sum [u_{j t}(z),u_{j t}(z)],
$$
where $(u_{j t})$ is an orthonormal basis for the scalar product
$H_t$. Explicitly
\be
B_t=\sum [u_j(z),u_j(z)]|t|^{-2\lambda_j}
\ee
if $(u_j)$ is  the diagonalizing basis chosen above.

For $p$ a positive integer we can do the same construction for
the space $\hat E(p)$ consisting of sections to $L^p\gr K_Z$ and get
metrics
$$
p\phi_{(p)}(t,\cdot)=\log B_t(p),
$$
on $L^p\gr K_Z$. 

Let $\phi_0$ and $\phi_1$ be two smooth positive metrics on $\hat L$, i e
points in $\K$. Let $U$ be the annulus $\{0<\log|t|<1\}$ and consider
the space $\mathcal{A}$ of all smooth semipositive metrics $\phi$ on $L$, the
pull back of $\hat L$ to 
$U\times Z$ such that
$$
\phi\leq \phi_0
$$
for $\log|t|=0$ and
$$
\phi\leq \phi_1
$$
for $\log|t|=1$. 

Then $\phi^*:=\sup_\A \phi$ is a moral candidate for a geodesic in $\K$,
but its eventual smoothness properties are a very hard issue, see
\cite{Chen}, \cite{Chen-Tian}.  

Here we shall prove a variant of the result of
Phong and Sturm, \cite{Phong-Sturm}. For $p$ large we consider the
Hilbert norms $H_{p\phi_0}$ and $H_{p\phi_1}$ on $\hat E(p)$, connect
them with a flat curve of Hilbert norms and define $B_t(p)$ and
$\phi_{(p)}$ as above. 
\begin{thm} 
\be
\sup |\phi_{(p)} - \phi^*|\leq C\frac{\log p}{p}
\ee
\end{thm}
The meaning of this statement is perhaps a bit obscure since $\phi_{(p)}$
and $\phi^*$ are metrics on different bundles. If we choose a fixed
smooth metric $\chi$ on $K_Z$ the precise meaning of (6.2) is
$$
\sup |(\phi_{(p)}-\chi/p) - \phi^*|\leq C\frac{\log p}{p}.
$$
Note that in the case of principal interest, when $K_Z$ is negative,
we can choose $\chi$ to have negative curvature, so that our
approximants are positively curved. In case $\hat L$ is a power of the
canonical bundle on $Z$ one can also avoid the introduction of $\chi$
by normalizing $\phi_{(p)}$ differently.

To prove Theorem 6.1 we first note that $\phi^*$ is bounded from above
by $\max(\phi_0,\phi_1)$. This follows if we apply the maximum
principle with resepct to the $t$-variable for a general element in $\A$ 

The direction of Theorem 6.1 that estimates $\phi^*$ from below is
relatively straightforward. First note that when $\log|t|=0$, by the
Tian-Zelditch-Catlin  formula
\be
|\phi_{(p)}-\chi/p -\phi_0|\leq C\log p/p
\ee
and that a similar estimate holds on the outer boundary of the
annulus. We will use $\phi_0$ as a fixed strictly positive auxilary
metric on $\hat L$. Notice that by (6.1) $\phi_{(p)}$ defines a
semipositive metric on $L$ over $U\times Z$, so that if $a$ is a
sufficiently large constant
$$
\phi_{(p)}(1-a/p) +a\phi_0/p -\chi/p
$$
is also semipositive for $p$ large enough (the positivity of $\phi_0$
compensates for the possible negativity of $\chi$ if $a$ is large
enough) . Combining with (6.3) we see 
that
$$
\xi_p =\phi_{(p)}(1-a/p) +a\phi_0/p -\chi/p -C\log p/p
$$
belongs to $\mathcal{A}$. Hence $\phi^*\geq\xi_p$, proving one direction of
(6.2). Notice that this shows in particular that $\phi^*$ is uniformly
bounded from below by the smooth metric $\xi_{p_0}$ for some fixed
large $p_0$.

The proof of the other direction is divided into two steps. First we
estimate $B_t(p)$ from below  by $B_{p\phi^*}$ - the Bergman kernel
associated to $p\phi^*$, and then we estimate $B_{p\phi^*}$ by
$e^{p\phi^*}$. For the first step we need a well known lemma, cf
\cite{Semmes 2}.
\begin{lma} Let $E$ be a holomorphic vector bundle over (the closure
  of) a
  one-dimensional domain $U$, and let $A$ and $B$ be two Hermitian
  metrics on $E$ that extend continuously to $\bar U$. Assume that the
  curvature of $A$ is seminegative and that the curvature of $B$ is
  semipositive. Then, if $A\leq B$ on the boundary of $U$, it follows
  that $A\leq B$ in $U$.
\end{lma}

Let $E=E(p)$ be our trivial bundle with fiber $\hat E(p)$ and let $A$ be the
metric defined by the flat curve between $H_{p\phi_0}$ and
$H_{p\phi_1}$. $B$ is the metric defined by $H_{\phi^*}$. By Theorem
2.1, $B$ is semipositive, and by definition $A$ is flat. Thus by the
lemma
$$
A\leq B.
$$
This implies the opposite inequality for the Bergman kernels, so
\be
B_t(p)\geq B_{\phi^*}.
\ee
One might object here that $\phi^*$ is not necessarily smooth, so
Theorem 2.1 can not be applied directly. This can be circumvented by
proving instead (6.4) with $\phi^*$ replaced by an
arbitrary element in $\A$, which will suffice for our purposes. 

The remaining part of the proof now follows from a variant of the
Ohsawa-Takegoshi extension theorem.
\begin{thm}
Let $L$ be a line bundle over a compact manifold, $Z$,  with a
positive metric $\phi_0$, and let $\phi$ be a semipositive metric on $L$
such that $\phi_0-\phi$ is uniformly bounded. Let $\chi$ be a fixed
smooth metric on $K_Z$. 

Then, for any point $x$ in $Z$ and any sufficiently large integer $p$,
there is a holomorphic section, $h$  to $L^p\gr K_Z$ such that
\be
|h(x)|^2\geq  e^{p\phi +\chi}(x)
\ee
and 
\be
\int [h,h]e^{-p\phi}\leq C,
\ee
where $C$ does not depend on $p$.  
\end{thm} 
Accepting this for a moment we first see how Theorem 6.1 follows. We
will apply Theorem 6.3 to $\phi^*$ (or to an arbitrary element in
$\A$) for one fixed value of $t$. We 
already know that $\phi^*$ is bounded from below by a smooth metric
and  from above by
$\max(\phi_0,\phi_1)$, so Theorem 6.3 does apply to the couple
$\phi^*$ and $\phi_0$. By the extremal characterization of Bergman
kernels it follows 
from Theorem 6.3 that
$$
B_{p\phi^*}\geq C e^{p\phi^* +\chi}.
$$
Combining with (6.4) we find
\be
B_t(p)\geq C e^{p\phi^* +\chi}
\ee
from where it follows  that    
$$
\phi_{(p)}\geq \phi^* +\chi/p +C/p.
$$  
so Theorem 6.1 follows. 

To prove Theorem 6.3, first choose a trivializing neighbourhood $W$  and
local coordinates centered at $x$. By the Ohsawa-Takegoshi extension
theorem for bounded domains in $\Cn$ we can find a section over $W$
satisfying (6.5) over $W$ with an integral estimate over $W$. A
standard argument, involving Hörmander $L^2$-estimates over $Z$ with
respect to a singular weight with a logarithmic pole at $x$ then
extends $h$ to a global section such that (6.5) still holds and 
$$
\int [h,h]e^{- (p-p_0)\phi+p_0\phi_0}\leq C.
$$
Since $\phi-\phi_0$ is assumed to be uniformly bounded this gives
(6.6) and we are done.

As a final remark we note that if we assume known that $\phi^*$ has a
certain degree of smoothness, and that $\omega^t>0$, then we can
replace the crude lower bound in (6.7) from the Ohsawa-Takegoshi
theorem, by a few terms from the Tian-Zelditch-Catlin expansion. One
then gets a very precise estimate from below of $B_t(p)$.

\bigskip

\def\listing#1#2#3{{\sc #1}:\ {\it #2}, \ #3.}

\end{document}